\documentclass{ifacconf}

\usepackage{graphicx}      % include this line if your document contains figures
\usepackage{natbib}        % required for bibliography
\usepackage{amssymb}
\usepackage{mathtools}
\usepackage{accents}
\newcommand{\ubar}[1]{\underaccent{\bar}{#1}}
\newcommand{\norm}[1]{\left\Vert#1\right\Vert} % Norm
\DeclareMathAlphabet{\mbf}{OT1}{ptm}{b}{n}
\DeclareMathOperator*{\argmax}{arg\,max}
\DeclareMathOperator*{\argmin}{arg\,min}

\begin{document}
\begin{frontmatter}

\title{Semi-Infinite Programs for Robust Control and Optimization: Efficient Solutions and Extensions to Existence Constraints} 
% Title, preferably not more than 10 words.

\author[EEE]{Jad Wehbeh} 
\author[EEE,Aero]{Eric C. Kerrigan}

\address[EEE]{Department of Electrical and Electronic Engineering, Imperial College London, UK (e-mail: j.wehbeh22@imperial.ac.uk, e.kerrigan@imperial.ac.uk)}
\address[Aero]{Department of Aeronautics, Imperial College London, UK}

\begin{abstract}              
Discrete-time robust optimal control problems generally take a min-max structure over continuous variable spaces, which can be difficult to solve in practice. In this paper, we extend the class of such problems that can be solved through a previously proposed local reduction method to consider those with existence constraints on the uncountable variables. We also consider the possibility of non-unique trajectories that satisfy equality and inequality constraints. Crucially, we show that the problems of interest can be cast into a standard semi-infinite program and demonstrate how to generate optimal uncertainty scenario sets in order to obtain numerical solutions. We also include examples on model predictive control for obstacle avoidance with logical conditions, control with input saturation affected by uncertainty, and optimal parameter estimation to highlight the need for the proposed extension. Our method solves each of the examples considered, producing violation-free and locally optimal solutions.
\end{abstract}


\end{frontmatter}
%===============================================================================

\section{Introduction}

Min-max optimal control problems are control problems that seek to achieve the best possible performance in a given objective function for the worst-case realization of some unknown quantities, while ensuring that constraints are satisfied for all possible values of the uncertainty. When these possible uncertainties are drawn from infinite (i.e., to be more precise, uncountable) sets, the constraints are required to hold at infinitely many points, which can be difficult to verify in practice. 

In \citet{blankenship1976infinitely}, the authors show how semi-infinite programs (SIPs), which are optimization problems with a finite number of decision variables and an infinite number of constraints, can be solved through a scenario-based approach. Their proposed algorithm, which we will refer to as a local reduction method, reduces the infinite-dimensional constraint space to an equivalent finite number of key scenarios, which is more tractable from a numerical point of view. 

More recently, the work of \citet{zagorowska2023automatic} demonstrated the use of local reduction methods as a means of solving robust optimal control problems. The method they suggested allowed for the solution of numerous problems that were previously considered prohibitively challenging, but their proposed formulation was incapable of dealing with issues such as existence constraints or non-uniqueness in the solutions to the variables involved. This paper demonstrates how a wider range of robust optimal control problems can be solved through local reduction.

\vspace{-4pt}
\subsection{State of the Art}
\vspace{-4pt}

Scenario-based approaches to the solution of robust control problems have seen a recent growth in popularity, as shown by the survey of \citet{campi2021scenario}. This is particularly true for the field of model predictive control (MPC), where scenario sampling allows for tractable modelling of disturbances and uncertainties, as seen in the early work of \citet{calafiore2012robust} on linear systems or the more current application of nonlinear MPC to  heating of buildings \citep{pippia2021scenario}.

Nonetheless, choosing which scenarios to include in the constraint set is far from obvious, with the simplest methods relying on random sampling from a chosen distribution, such as with the uniformly distributed scenarios used in obtaining probabilistic guarantees by \citet{grammatico2015scenario}. Other approaches choose instead to generate scenarios by using scenario trees, which usually consider the lower bound, upper bound, and midpoint of each uncertain quantity at every discretization time step. This can be seen in the work of \citet{thangavel2018dual} on nonlinear MPC for dual control, or in  \citet{thombre2021sensitivity}, which further simplifies the problem by categorizing scenarios into critical or non-critical ones, and replacing the non-critical constraints by a sensitivity-based approximation.

However, worst-case realizations of the uncertainty are not guaranteed to lie on the edge or center of the uncertainty sets, as demonstrated in \citet{puschke2018robust}. This implies that scenario trees do not necessarily cover the worst-case uncertainty realizations, and may be required to include a prohibitively large number of scenarios to closely approximate this worst case. This issue is addressed by local reduction methods such as those of \citet{hettich1993semi}, which solve a separate maximization problem to find the worst case violations and add them to the constraint set, a fact that is exploited by \citet{zagorowska2023automatic} to effectively approximate the infinite-dimensional constraints on robust optimal control problems. 

Of course, it is worth noting that global optimization solvers can also be used to solve SIPs, as shown by \citet{djelassi2021recent}. These solvers guarantee that the solution obtained is actually a global minimum, an assertion that cannot be made when local solvers are used as part of local reduction methods, but suffer from the curse of dimensionality, which renders them unusable for larger control problems over longer prediction horizons due to their prohibitive computational cost.

\vspace{-4pt}
\subsection{Contributions}
\vspace{-4pt}

In this paper, we demonstrate how a wide range of robust optimal control problems can be cast as semi-infinite programs of the form considered in \citet{blankenship1976infinitely}. Specifically, we extend the work of \citet{zagorowska2023automatic} to cover problems with:
\begin{itemize}
    \item Existence constraints on non-unique sets of variables defining the constraints, which can arise from the smoothing of logical conditions and transform the problem into a min-max-min.
    \item Non-uniqueness in the variables used for the propagation of the state, which can arise as a result of the introduction of additional modelling variables that facilitate the expression of the dynamics. 
    \item Systems of equalities and inequalities of the state, control, and input variables governing the feasible evolution of the state trajectories. 
\end{itemize}

The paper begins by introducing motivational examples for each of the proposed contributions in Section \ref{sec:exs} before detailing the problem formulation in Section \ref{sec:form}. Subsequently, an overview of how local reduction can be used to solve the problem considered is included in Section \ref{sec:local_red}, with numerical results shown in Section \ref{sec:results}.

\section{Motivational Examples}
\label{sec:exs}

In order to illustrate the motivation behind some of the extensions being suggested, we now consider three different examples focusing on separate problem aspects. The examples are intentionally chosen to be simple so as to better highlight the parts of the problem definition that can be an issue for existing solution approaches. 

\vspace{-4pt}
\subsection{Obstacle Avoidance}
\label{sec:ex1}
\vspace{-4pt}

To begin with, we introduce a basic trajectory planning problem for a three-dimensional system under unknown but bounded disturbances. The system of interest has the discrete-time linear dynamics 
\begin{equation}
    x_{i,k+1} = x_{i,k} + u_{i,k} + w_{i,k} \qquad \forall i \in \{1,2,3\}
\end{equation}
where $x_{i,k}$ is the $i$th state at time step $k$, $u_{i,k}$ is the $i$th control input at time step $k$, and $w_{i,k}$ is the $i$th bounded disturbance at time step $k$.

Given an initial state $\mbf{x}_0$, the goal of the problem is to navigate to a target position $\mbf{x}_{ref}$, aiming for the best possible performance under the worst-case realization of the uncertainty. The objective also penalizes the total control input applied, as described by
\begin{equation}
\label{eq:ex1_obj}
    \min_{\mbf{u}} \max_{\mbf{x},\mbf{w}} \left( \; \sum_{k=0}^{n-1} \norm{\mbf{u}_k}^2_R + \norm{\mbf{x}_n-\mbf{x}_{ref}}^2_Q \right).
\end{equation}
Here, bolded symbols are used to describe vectors of the bolded quantity correspondoing to a specific time step. For instance, $\mbf{x}_k := \begin{bmatrix} x_{1,k} & x_{2,k} & x_{3,k} \end{bmatrix}^\top$. Also, $\norm{\mbf{x}}_Q^2 := \mbf{x}^\top Q \,\mbf{x}$.

%TODO: Add figure

On top of minimizing the objective of \eqref{eq:ex1_obj}, the system is also required to avoid a cylindrical obstacle by satisfying the constraint
\begin{equation}
    \label{eq:ex1_constr_nonsmooth}
    x_{1,k}^2 + x_{2,k}^2 \geq 1 \quad \lor \quad x_{3,k} \geq 1 \quad \lor \quad x_{3,k} \leq -1 
\end{equation}
at every time step $k$, where $\lor$ is the logical OR operator. This is equivalent to stating that the system must always lie outside, above, or below the cylindrical obstacle.

Of course, giving a constraint of the form of \eqref{eq:ex1_constr_nonsmooth} directly to a solver can be rather problematic, so we choose to reformulate it as
\begin{equation}
    \min \left(1 - x_{1,k}^2 + x_{2,k}^2,\, 1 - x_{3,k},\, 1 + x_{3,k}\right) \leq 0
\end{equation}
which ensures that any of the three conditions is satisfied. We then apply the smoothing approach of \citet{kirjner1998conversion}, resulting in the equivalent smoothed constraint
\begin{equation}
    \label{eq:ex1_constr_smooth}
    s_{1,k}(1 - x_{1,k}^2 + x_{2,k}^2) + s_{2,k} (1 - x_{3,k}) + s_{3,k}(x_{3,k} +1) \leq 0 
\end{equation}
where the smoothing variables must satisfy
\begin{subequations}
    \label{eq:eq:ex1_constr_cond}
    \begin{align}
        s_{1,k} + s_{2,k} + s_{3,k} = 1 \\
        s_{1,k}, s_{2,k}, s_{3,k} \geq 0.
    \end{align}
\end{subequations}

For \eqref{eq:ex1_constr_smooth} to guarantee that the constraint of \eqref{eq:ex1_constr_nonsmooth} is satisfied, it is simply necessary and sufficient for there to exist a choice of $\mbf{s}_k$ satisfying \eqref{eq:ex1_constr_smooth} and \eqref{eq:eq:ex1_constr_cond} for every possible realization of the uncertainty. However, it is trivial to observe that while \eqref{eq:ex1_constr_smooth} and \eqref{eq:eq:ex1_constr_cond} can admit non-unique solutions $\mbf{s}_k$, it is only necessary to prove the existence of a single solution for the original constraint of \eqref{eq:ex1_constr_nonsmooth} to be satisfied. This is a problem that cannot be solved using existing local reduction-based SIP solvers.

\vspace{-4pt}
\subsection{Input Saturation}
\label{sec:ex2}
\vspace{-4pt}

Next, we examine the exact modelling of saturation constraints on the control input. In order to do so, we consider the system dynamics
\begin{equation}
    \label{eq:ex2_dyn}
    x_{k + 1} = (a + w)x_k + u_k
\end{equation}
where $a > 1$ is a constant characterizing the unstable dynamics of $x$, and $w$ is some bounded uncertainty on the true value of $a$. The control input is computed through the feedback law $u_k = - b x_k$ where $b$ is a fixed gain, and saturated to be between $\ubar{u}$ and $\bar{u}$.

In the previous work of \citet{zagorowska2023automatic}, input saturation was handled approximately through the use of an exponential saturation function. Here, we choose to instead consider the explicit saturation
\begin{equation}
\label{eq:sat_explicit}
    u_k = \text{sat}(\ubar{u},\bar{u},-b x_k)
\end{equation}
where $\text{sat}(a_1,a_2,a_3)$ returns the projection of $a_3$ in the interval $[a_1, a_2]$ with $a_1 \leq a_2$. Note that the saturation constraint of \eqref{eq:sat_explicit} cannot simply be replaced by $\ubar{u} \leq -b x_k \leq \bar{u}$, because doing so would limit the space of admissible values for $b$ instead of saturating $u_k$.

\begin{prop}
\label{prop:saturation}
The constraint of \eqref{eq:sat_explicit} can be substituted exactly for the smoothed formulation
\begin{subequations}
\label{eq:sat_smooth}
    \begin{align}
        z_{1,k}(-b x_k - \bar{u}) + z_{2,k}(u_k - \bar{u})^2 & \leq 0 \\
        z_{3,k}(b x_k + \ubar{u}) + z_{4,k}(u_k - \ubar{u})^2 & \leq 0 \\
        \begin{split}
            z_{5,k}(b x_k + \bar{u}) + z_{6,k}(-b x_k - \ubar{u}) \; &  \\
                + z_{7,k}(u_k + b x_k)^2 & \leq 0
        \end{split}
    \end{align}
where each $z_{i,k}$ is non-negative and
    \begin{align}
        z_{1,k} + z_{2,k} = 1 \\
        z_{3,k} + z_{4,k} = 1 \\
        z_{5,k} + z_{6,k} + z_{7,k} = 1
    \end{align}
\end{subequations}
\end{prop}

\begin{pf}
Given the premises and conclusions
\begin{equation*}
    \begin{array}{ll}
     p_1 : -b x_x \leq \bar{u}  & \qquad c_1 : u_k = - b x_k \\
     p_2 : -b x_x \geq \ubar{u} & \qquad c_2 : u_k = \bar{u} \\
                                & \qquad c_3 : u_k = \ubar{u}
    \end{array}
\end{equation*}
the constraint of \eqref{eq:sat_explicit} is equivalent to the logical statement
\begin{equation}
\label{eq:sat_constr_log_impl}
    (\neg \, p_1 \implies c_2) \land (\neg \, p_2 \implies c_3) \land [(p_1 \land p_2) \implies c_1]   
\end{equation}
which can be transformed using the equivalence relation $(p \implies c) \iff (\neg \, p \lor c)$ to
\begin{equation}
\label{eq:sat_constr_log_or}
    (p_1 \lor c_2) \land (p_2 \lor c_3) \land [\neg \, p_1 \lor \neg \, p_2 \lor c_1]
\end{equation}

The first condition, $(p_1 \lor c_2)$, can be expanded into
\begin{equation}
    ( -b x_x \leq \bar{u} ) \lor ( u_k = \bar{u})
\end{equation}
and reformulated using minimums into
\begin{equation}
\label{eq:sat_const_min_1}
    \min\left\{-b x_k - \bar{u}, (u_k - \bar{u})^2\right\} \leq 0.
\end{equation}

Using the smoothing approach of \citet{kirjner1998conversion} once more, the constraint of \eqref{eq:sat_const_min_1} becomes
\begin{subequations}
\label{eq:sat_const_smooth_1}
\begin{equation}
    z_{1,k}(-b x_k - \bar{u}) + z_{2,k}(u_k - \bar{u})^2  \leq 0
\end{equation}
with
\begin{align}
    z_{1,k} + z_{2,k} & = 1 \\
    z_{1,k}, \, z_{2,k} & \geq 0.
\end{align}
\end{subequations}

The same procedure can then be applied to the two other conditions in \eqref{eq:sat_constr_log_or} to recover \eqref{eq:sat_smooth}. Since every step between \eqref{eq:sat_constr_log_impl} and \eqref{eq:sat_const_smooth_1} is a two-way implication, \eqref{eq:sat_explicit} and \eqref{eq:sat_smooth} must be equivalent, completing the proof. \qed
\end{pf}

The goal of this example is to find a $b$ that solves
\begin{equation}
    \min_{b} \max_{\mbf{u}, \, \mbf{x}, \, \mbf{z}, \, w} \norm{x_n}^2
\end{equation}
subject to \eqref{eq:ex2_dyn} and \eqref{eq:sat_smooth} at some time step $n > 0$ given a fixed initial condition $x_0$ and a bounded $w$. Once again, any of a non-unique set of possible solutions for $\mbf{z}_k$ must be found in order to propagate the dynamics, but this cannot be handled by existing robust optimal control methods. 

\vspace{-4pt}
\subsection{Parameter Estimation}
\label{sec:ex3}
\vspace{-4pt}

Finally, we pose the problem of robust mass estimation for the discretized double integrator described by  
\begin{subequations}
\label{eq:ex3_dyn}
    \begin{align}
        x_{1,k+1} & = x_{1,k} + \Delta t \, x_{2,k} \\
        \label{eq:ex3_dyn_x2}
        x_{2,k+1} & = x_{2,k} + \Delta t \, \frac{u_{k}}{m} 
    \end{align} 
\end{subequations}
where $m$ represents the mass of the system, $x_{1,k}$ is the position at time step $k$, $x_{2,k}$ is the velocity, and $u_k$ is the force applied. $\Delta t$ is the discretization time step used.

Measurements $y_k$ of $x_{1,k}$, which are corrupted by a noise~$w_k$, are taken at each time step and used directly in obtaining the estimate $\hat{x}_{1,k}$ through
\begin{equation}
\label{eq:ex3_est_x1}
    \hat{x}_{1,k} = y_k = x_{1,k} + w_k.
\end{equation}       

The estimate $\hat{x}_{2,k}$ of $x_{2,k}$ is then simply obtained by taking a finite difference between successive measurements, resulting in the approximation
\begin{equation}
\label{eq:ex3_est_x2}
\begin{split}
     \hat{x}_{2,k} & = \left(y_{k+1} - y_k\right)/\Delta t = \left(\hat{x}_{1,k+1} - \hat{x}_{1,k}\right)/\Delta t \\ & = \left({x}_{1,k+1} - {x}_{1,k} + {w}_{1,k+1} - {w}_{1,k}\right)/\Delta t \\
    & = x_{2,k} + \left({w}_{1,k+1} - {w}_{1,k}\right)/\Delta t.
\end{split}
\end{equation}

In this parameter estimation problem, we seek to obtain a lower bound $\underline{m}$ and an upper bound $\overline{m}$ for the unknown value of $m$, such that
\begin{equation}
    \label{eq:ex3_constr_m}
    \underline{m} \leq m \leq \overline{m}.
\end{equation}

Using \eqref{eq:ex3_dyn_x2} and $\eqref{eq:ex3_constr_m}$, we know $\overline{m}$ must satisfy
\begin{subequations}
\label{eq:ex3_constr_mbar}
    \begin{align}
        \frac{\Delta t}{\overline{m}} & \leq \frac{{x}_{2,k+1} - {x}_{2,k}}{u_k}\\
        \frac{\Delta t}{\overline{m}} & \leq \frac{\hat{{x}}_{2,k+1} - \hat{{x}}_{2,k+1} - \frac{w_{k+2} - 2w_{k+1} + w_k}{\Delta t}}{u_k} \\
        \frac{\Delta t^2}{\overline{m}} & \leq \frac{y_{k+2} - 2y_{k+1} + y_{k} - w_{k+2} + 2w_{k+1} - w_k}{u_k}.
    \end{align}
\end{subequations}

Similarly, $\underline{m}$ must satisfy
\begin{equation}
\label{eq:ex3_constr_mubar}
    \frac{\Delta t^2}{\underline{m}} \geq \frac{y_{k+2} - 2y_{k+1} + y_{k} - w_{k+2} + 2w_{k+1} - w_k}{u_k}.
\end{equation} 

In order to obtain the tightest possible bounds on the estimate of $m$, and given a series of measurements $\mbf{y}$ and known control inputs $\mbf{u}$, we target the objective
\begin{equation}
    \min_{\underline{m},\overline{m}} \max_{\mbf{x}, \, \mbf{w}} \; \overline{m} - \underline{m}
\end{equation}
subject to the constraints of \eqref{eq:ex3_dyn}, \eqref{eq:ex3_est_x1}, \eqref{eq:ex3_constr_mbar}, and \eqref{eq:ex3_constr_mubar}.

At first glance, this parameter estimation example appears to be a standard min-max problem of the type considered by \citet{zagorowska2023automatic}. However, in this problem, the equality and inequality constraints combine to limit the feasible combinations of states and uncertainties, which \citet{zagorowska2023automatic} does not consider.

\section{Problem Formulation}
\label{sec:form}

\vspace{-4pt}
\subsection{System Definition}
\label{sec:form_sys_def}
\vspace{-4pt}

Consider the discrete-time system with state trajectory $z$ $\in$ $\mathcal{Z}$ $\subseteq$ $\mathbb{R}^{n_z}$, which can be split into physical states $z_p$ $\in$ $\mathcal{Z}_p$ $\subseteq$ $\mathbb{R}^{n_p}$ and modelling variables $z_m$ $\in$ $\mathcal{Z}_m$ $\subseteq$ $\mathbb{R}^{n_m}$, such that $n_z = n_m + n_p$. Here, $z_m$ represents variables that are introduced for the purpose of facilitating the modelling of the system but are not necessarily bound by physical laws, such as the $z$ variables in the example of Section \ref{sec:ex2}.

The evolution of $z$ is governed by the the equations $d(\theta,w,z_p,z_m) = 0$, which may arise from the explicit or implicit discretization of the continuous-time system dynamics, and where $\theta$ $\in$ $\Theta$ $\subseteq$ $\mathbb{R}^{n_\theta}$ is a collection of decision variables including control inputs and gains, and $w$ $\in$ $\mathcal{W}$ $\subseteq$ $\mathbb{R}^{n_w}$ represents the uncertainties the system is subject to. The natural evolution of the trajectories is also limited by a set of inequalities $e(\theta,w,z_p,z_m) \leq 0$ that are required to hold at all times. The equations involving $d$ and $e$ allow the choices of $\theta$ and $w$ to constrain the evolution of the state trajectory, enabling the transcription of the examples of Sections \ref{sec:ex2} and \ref{sec:ex3}. The trajectory space $\mathcal{Z}$ is therefore  represented as 
\begin{equation}
\mathcal{Z}(\theta,w) \coloneqq \left\{ (z_p,z_m) \left| \:
    \begin{array}{l}
        d(\theta,w,z_p,z_m) = 0 \\
        e(\theta,w,z_p,z_m) \leq 0 
    \end{array}
    \right. \right\}
\end{equation}
where $\Theta$ and $\mathcal{W}$ are compact and defined independently. Note that in the case where the discrete-time system considered is obtained from the discretization of a continuous-time problem, the values of $n_\theta$, $n_w$, and $n_z$ depend on all of the dimensions of the respective continuous variables, the number of time steps considered, and the discretization method used. 

Given the intended natures of $z_p$ and $z_m$, we introduce the following assumption:

\begin{assum}
\label{assum:uniqueness}
$\forall$ $(\theta,w)$ $\in$ $\Theta \times \mathcal{W}$, there exists a unique $z_p$ and a set $M$ such that $\forall$ $z_m$ $\in$ $M$: $d(\theta,w,z_p,z_m) = 0$ and $e(\theta,w,z_p,z_m) \leq 0$. Additionally, $M \neq \varnothing$ $\forall$ $(\theta,w)$ $\in$ $\Theta \times \mathcal{W}$.  
\end{assum}

Assumption \ref{assum:uniqueness} is guaranteed to be justified when $z_p$ is used to describe components of the state that evolve along unique trajectories in continuous-time, and the discretization method used preserves this uniqueness.
 
Under this assumption, we can define the functions $z_p^*(\cdot,\cdot):$ $\Theta \times \mathcal{W} \rightarrow \mathcal{Z}_p$ and $z_m^*(\cdot,\cdot):$ $\Theta \times \mathcal{W} \rightarrow 2^{\, \mathcal{Z}_m}$, which map any choice of $\theta$ and $w$ onto the corresponding values of $z_p$ and $M$, and where $2^{\, \mathcal{Z}_m}$ is the power set of $\mathcal{Z}_m$. Furthermore, under Assumption \ref{assum:uniqueness}, we can show that $\mathcal{Z}(\theta,w) = \mathcal{Z}_p(\theta,w) \times \mathcal{Z}_m(\theta,w)$, where 
\begin{align}
\begin{split}
&\mathcal{Z}_p(\theta,w) \coloneqq \left\{ z_p \left| \, z_p = z_p^*(\theta,w) \right. \right\} \coloneqq \\
&\quad \left\{ z_p \left| \, \exists z_m \in z^*_m(\theta,w) : \:
    \begin{array}{l}
        d(\theta,w,z_p,z_m) = 0 \\
        e(\theta,w,z_p,z_m) \leq 0 
    \end{array}
    \right. \right\}
\end{split}\\[5pt]
\begin{split}
& \mathcal{Z}_m(\theta,w) \coloneqq \left\{ z_m \left| \, z_m \in M, \, M = z_m^*(\theta,w) \right. \right\} \coloneqq \\
&\quad \left\{ z_m \left| \:
    \begin{array}{l}
        d(\theta,w,z_p^*(\theta,w),z_m) = 0 \\
        e(\theta,w,z_p^*(\theta,w),z_m) \leq 0 
    \end{array}
    \right. \right\} .
\end{split}
\end{align}

We also make the assumption:

\begin{assum}
\label{assum:zp_cont}
The unique mapping for $z_p^*$ is continuous, i.e., $z_p^*(\cdot,\cdot) \in \mathcal{C}^{\,0}$.
\end{assum}

\vspace{-4pt}
\subsection{Robust Optimal Control Problem}
\label{sec:form_contr_prob}
\vspace{-4pt}

We now define a robust optimal control problem for the system of Section \ref{sec:form_sys_def} that generalizes the formulation of \citet{zagorowska2023automatic}. To do so, we first consider the inner problem
\begin{subequations}
\begin{equation}
    f^*(\theta) \coloneqq \max_{w \in \mathcal{W}, \, z_p \in \mathcal{Z}_p(\theta,w)} f(\theta,w,z_p)
\end{equation}
\end{subequations}
which finds the worst-case value of the objective function $f$ under any uncertainty realization $w$ and associated feasible state trajectory $z_p$. A robust optimal control problem for the system can then be formulated by minimizing $f^*(\theta)$ over the admissible values of $\theta$, while ensuring a vector of constraints $g$ is satisfied for all realizations of the uncertainty. This can be described by the problem
\begin{subequations}
\label{eq:rob_prob}
    \begin{equation}
        \min_{\theta \in \Theta} f^*(\theta)
    \end{equation}
    s.t. $\forall w \in \mathcal{W}, \forall z_p \in \mathcal{Z}_p(\theta,w)$
    \begin{equation}
        \label{eq:s_g_constr}
        \exists s \in S: g(\theta,w,z_p,s) \leq 0 
    \end{equation}
\end{subequations}
where
\begin{equation}
\mathcal{S} \coloneqq \left\{ s \left| \:
    \begin{array}{l}
        q(s) = 0 \\
        r(s) \leq 0 
    \end{array}
    \right. \right\}
\end{equation}
is the space of admissible values for $s$, a variable for which there must exist a solution to ensure that \eqref{eq:s_g_constr} is satisfied. This constraint formulation is motivated by the example of Section \ref{sec:ex1}, for which a value of $s$ satisfying the constraints must exist in order to ensure the obstacles are successfully avoided. Note that $s$ differs from $z_m$ in that $s$ must exist for the inequality constraints to be satisfied while $z_m$ must exist for the propagation of the dynamics to be valid.

We additionally assume:
\begin{assum}
\label{assum:s_compact}
The set $\mathcal{S}$ is compact.
\end{assum}
\begin{assum}
\label{assum:fg_cont}
The functions $f$ and $g$ are continuous in all of their arguments, i.e., $f(\cdot,\cdot,\cdot), \, g(\cdot,\cdot,\cdot,\cdot) \in \mathcal{C}^{\,0}$.
\end{assum}
\begin{assum}
\label{assum:obj_bounded}
The value of $f(\theta,w,z_p^*(\theta,w))$ is bounded $\forall \theta \in \Theta$, $\forall w \in \mathcal{W}$.
\end{assum}
\begin{assum}
\label{assum:constr_bounded}
The value of $g(\theta,w,z_p^*(\theta,w),s)$ is bounded $\forall \theta \in \Theta$, $\forall w \in \mathcal{W}$, $\forall s \in \mathcal{S}$.
\end{assum}

Using Assumption \ref{assum:uniqueness}, \eqref{eq:rob_prob} can be rewritten as
\begin{subequations}
    \begin{equation}
        \min_{\theta \in \Theta} f^*(\theta)
    \end{equation}
    s.t. $\forall w \in \mathcal{W}$
    \begin{equation}
        \exists s \in S: g(\theta,w,z_p^*(\theta,w),s) \leq 0.
    \end{equation}
\end{subequations}

Converting the problem to its subgraph form, this is equivalent to 
\begin{subequations}
\label{eq:robust_prob_full}
    \begin{equation}
        \min_{\theta \in \Theta, \, \gamma \in \Gamma} \gamma
    \end{equation}
    s.t. $\forall w \in \mathcal{W}$
    \begin{align}
        \label{eq:robust_prob_full_sconst}
        \exists s \in S: g(\theta,w,z_p^*(\theta,w),s) & \leq 0 \\
        f(\theta,w,z_p^*(\theta,w)) - \gamma & \leq 0
    \end{align}
\end{subequations}
where $\gamma \in \Gamma \subset \mathbb{R}$ is introduced as an upper bound on the objective value. \citet{djelassi2021global} introduces a method to deal with the existence constraint of \eqref{eq:robust_prob_full_sconst} explicitly, but relies on global optimization techniques and is unsuitable for adaptation to local reduction. Instead, we propose rewriting \eqref{eq:robust_prob_full} as
\begin{subequations}
\label{eq:robust_prob_min}
    \begin{equation}
        \min_{\theta \in \Theta, \, \gamma \in \Gamma} \gamma
    \end{equation}
    s.t. $\forall w \in \mathcal{W}$
    \begin{align}
        \label{eq:robust_prob_min_gconst}
        \min_{s \in S} \; g(\theta,w,z_p^*(\theta,w),s) & \leq 0 \\
        \label{eq:robust_prob_min_fconst}
        f(\theta,w,z_p^*(\theta,w)) - \gamma & \leq 0
    \end{align}
\end{subequations}
where \eqref{eq:robust_prob_min_gconst} is logically equivalent to \eqref{eq:robust_prob_full_sconst}.

From Assumptions \ref{assum:zp_cont} and \ref{assum:fg_cont}, we know that $(\theta, w, \gamma) \mapsto f(\theta, w, z_p^*(\theta, w)) - \gamma$ is continuous. Similarly, using Assumption \ref{assum:s_compact} for the compactness of $\mathcal{S}$, it follows from \citet{clarke1975generalized} that $(\theta, w) \mapsto \min_{s \in \mathcal{S}} g(\theta, w, z_p^*(\theta, w),s)$ is also continuous. Finally, given that $\Theta$ and $\mathcal{W}$ are compact by definition and $\Gamma$ is compact under Assumption \ref{assum:obj_bounded}, the problem of \eqref{eq:robust_prob_min} can be rewritten as 
\begin{equation}
    \min_{\bar{\theta}\in\bar{\Theta}} \bar{f}(\bar{\theta}) \quad \text{s.t.} \; \forall w \in \mathcal{W}: \; \bar{g}(\bar{\theta},w) \leq 0
\end{equation}
which is an SIP satisfying all of the assumptions of \citet{blankenship1976infinitely}. Here, $\bar{\theta}$ combines $\theta$ and $\gamma$, $\bar{f}(\bar{\theta}) = \gamma$, and $\bar{g}(\cdot,\cdot)$ collects \eqref{eq:robust_prob_min_gconst} and \eqref{eq:robust_prob_min_fconst}. 

\section{Local Reduction Algorithm}
\label{sec:local_red}
While the problem of \eqref{eq:robust_prob_min} is useful in proving that local reduction can be employed to solve the desired min-max-min problem, it is far from ideal for use with numerical solvers. Notably, the function $z_p^*$ is difficult to find explicitly, rendering the exact representation of \eqref{eq:robust_prob_min} close to impossible. Instead, we apply the local reduction algorithm described in \citet{zagorowska2023automatic} optimizing for the tuple $(\theta, \gamma)$ under a finite number of scenarios $w_i \in \mathbb{W}_j$ where $\text{card}\left(\mathbb{W}_j\right) = j$. Under any given policy $(\theta, \gamma)$, the worst-case uncertainty $w$ for the constraint of \eqref{eq:robust_prob_min_fconst} is found through
\begin{equation}
\label{eq:max_prob_fconst}
   \argmax_{w \in \mathcal{W}, \, z_P \in \mathcal{Z}_p(\theta,w)} f(\theta,w,z_p) - \gamma
\end{equation}
while that for each $i$-th row of \eqref{eq:robust_prob_min_gconst} is obtained from
\begin{equation}
\label{eq:max_prob_gconst_full}
   \argmax_{w \in \mathcal{W}, \, z_P \in \mathcal{Z}_p(\theta,w)} \min_{s\in\mathcal{S}} g_i(\theta,w,z_p,s).
\end{equation}

The problem of \eqref{eq:max_prob_gconst_full} can be substituted for
\begin{subequations}
\label{eq:max_prob_gconst_sub}
    \begin{equation}
       \argmax_{w \in \mathcal{W}, \, z_P \in \mathcal{Z}_p(\theta,w), \, \sigma} \sigma
    \end{equation}
s.t. $\forall s$ $\in$ $\mathcal{S}$
    \begin{equation}
        \sigma - g_i(\theta,w,z_p,s) \leq 0
    \end{equation}
\end{subequations}
which is a nested SIP that can be solved directly through local reduction under Assumption \ref{assum:constr_bounded}, since \eqref{eq:max_prob_gconst_sub} does not have any inner optimization problems inside its constraints in the way that \eqref{eq:robust_prob_min} does. 

Once the solutions of \eqref{eq:max_prob_fconst} and \eqref{eq:max_prob_gconst_sub} are obtained, a subset of the violation values of $w$ is added to $\mathbb{W}_j$, resulting in the new scenario set $\mathbb{W}_k$ such that $k > j$. The values of $z_p$ obtained as part of the solution are simply discarded, because  they are unnecessary in computing the new values of $(\theta, \gamma)$ that are obtained by solving 
\begin{subequations}
\label{eq:lr_min_prob}
    \begin{equation}
        \argmin_{\theta \in \Theta, \, \gamma} \gamma 
    \end{equation}
s.t. $\forall i$ $\in$ $\{1,\hdots,k\}$, $\exists\, z_{p,i},\,z_{m,i}, \, s_i:$
    \begin{align}
        \label{eq:lr_min_prob_g}
        g(\theta,w_i,z_{p,i},s_i) & \leq 0 \\
        \label{eq:lr_min_prob_f}
        f(\theta,w_i,z_{p,i}) - \gamma & \leq 0 \\
        \label{eq:lr_min_prob_d}
        d(\theta,w_i,z_{p,i},z_{m,i}) & = 0 \\
        e(\theta,w_i,z_{p,i},z_{m,i}) & \leq 0 \\
        q(s_i) & = 0 \\
        \label{eq:lr_min_prob_r}
        r(s_i) & \leq 0
    \end{align}
\end{subequations}
which is equivalent to
\begin{subequations}
\label{eq:lr_min_min_prob}
    \begin{equation}
        \argmin_{\theta \in \Theta, \, \gamma} \gamma 
    \end{equation}
s.t. $\forall i$ $\in$ $\{1,\hdots,k\}$,
    \begin{equation}
        \label{eq:lr_min_min_prob_g}
        \min_{z_{p,i},\,z_{m,i},\,s_i}  g(\theta,w_i,z_{p,i},s_i) \leq 0
    \end{equation}
\end{subequations}
along with the constraints \eqref{eq:lr_min_prob_f} through \eqref{eq:lr_min_prob_r}. The $\min$ in the constraint of \eqref{eq:lr_min_min_prob} can then be moved into the objective to reformulate Problem \eqref{eq:lr_min_min_prob_g} as
\begin{equation}
    \label{eq:lr_min_prob_practice}
    \argmin_{\theta \in \Theta, \, \gamma} \min_{z_{p,i},\,z_{m,i},\,s_i \; \forall i \in \{1,\hdots,k\}} \gamma 
\end{equation}
with the constraints of \eqref{eq:lr_min_prob_g} through \eqref{eq:lr_min_prob_r} holding $\forall i$ $\in$ $\{1,\hdots,k\}$. This is a finite-dimensional min problem that can be fed directly to a numerical solver.

The algorithm then iterates between adding new scenarios and recomputing $\theta$ and $\gamma$ until no values of $w$ resulting in constraint violations can be found by solving \eqref{eq:max_prob_fconst} and \eqref{eq:max_prob_gconst_sub}, implying that the most recent point for $(\theta,\gamma)$ is optimal. Alternatively, if the number of scenarios considered surpasses some allowable threshold, the algorithm can simply return a suboptimal choice of $(\theta,\gamma)$ that is specifically robust to all of the scenarios in $\mathbb{W}_k$.

Under this formulation, the values of $z_p$ and $z_m$ obtained in solving the problems of \eqref{eq:max_prob_fconst}, \eqref{eq:max_prob_gconst_sub}, and \eqref{eq:lr_min_prob} are completely independent, eliminating any issues associated with non-uniqueness. The problems simply require that feasible values of the state variables exist so as to ensure that the propagation is valid.

\section{Results}
\label{sec:results}

We now provide the solutions obtained for the three examples described in Section \ref{sec:exs} after applying the algorithm of Section \ref{sec:local_red}. The problems and algorithm are implemented in the Julia programming language using the JuMP package \citep{Lubin2023} and the Ipopt optimizer \citep{wachter2006implementation}. Since Ipopt is a local solver, however, the solutions obtained to \eqref{eq:max_prob_fconst}, \eqref{eq:max_prob_gconst_sub}, and \eqref{eq:lr_min_prob_practice} cannot be guaranteed to be globally optimal, eliminating any guarantees on the global optimality of the overall solutions obtained. Instead, we rely on random restarts to ensure the eventual discovery of all relevant uncertainty scenarios and the convergence to a feasible solution. This enables ``good enough" solutions to be obtained much faster than through global optimization approaches. Each example is initialized with a single scenario, for which all of the uncertainties are chosen to be zero.

\vspace{-4pt}
\subsection{Obstacle Avoidance}
\vspace{-4pt}

For the obstacle avoidance example of Section \ref{sec:ex1}, the initial position is chosen to be $\mbf{x}_0 = \begin{bmatrix} -2 & 0 & 0 \end{bmatrix}^\top$ and the target position is $\mbf{x}_f = \begin{bmatrix} 2 & 0 & 0 \end{bmatrix}^\top$ after $n = 5$ time steps. We also choose $\mbf{Q} = \mbf{1}_{3 \times 3}$ and $\mbf{R} = 0.05\times\mbf{1}_{4 \times 4}$, where $\mbf{1}_{m \times m}$ is the identity matrix of size $m$. The control inputs and disturbances are limited to the range
\begin{subequations}
\begin{align}
    -1.0 \leq & u_{i,k} \leq 1.0 && \forall i \in \{1,2,3\} \\
    -0.1 \leq & w_{i,k} \leq 0.1 && \forall i \in \{1,2,3\}.
\end{align}
\end{subequations}

\begin{figure}[b]
    \centering
    \includegraphics[width = 0.72\columnwidth]{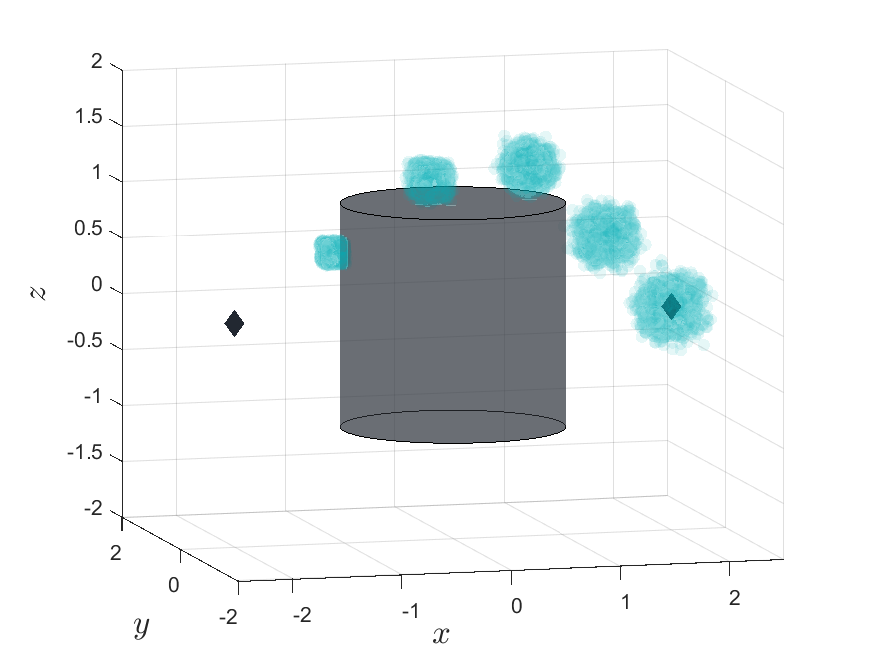}
    \caption{Distribution of 1000 randomly sampled trajectories under the open loop control inputs obtained}
    \label{fig:ex1_plot}
\end{figure}

The optimal open-loop control trajectory obtained after the addition of 29 scenarios is
\begin{equation}
    \mbf{u}^* = \begin{bmatrix}
        0.9 & 0.9 & 0.9 & 0.7 & 0.6 \\ 0.0 & 0.0 & 0.0 & 0.0 & 0.0 \\ 0.6 & 0.6 & 0.1 & -0.65 & -0.65
    \end{bmatrix}
\end{equation}
which successfully avoids collisions with the obstacle when simulated under uniformly randomly sampled values for $w_{i,k}$, as shown in  Figure \ref{fig:ex1_plot}. The solution procedure also produces the upper bound on the cost $\gamma^{\, *} = 0.714$, which is larger than the largest cost of 0.673 produced after one million random runs of the simulation.

\vspace{-4pt}
\subsection{Input Saturation}
\vspace{-4pt}

The second example we solve is that of Section \ref{sec:ex2} with $n = 5$ time steps. We take $a = 1.3$ and constrain $w$ to be between -0.2 and 0.2. Running the  algorithm results in convergence after the addition of a single scenario for which $w = 0.2$. The optimal control gain is determined to be $b^* = 1.339$ and the upper bound on the cost obtained is $\gamma^{\, *} = 8.94 \times 10^{-6}$.

Simulating the dynamics of $\eqref{eq:ex2_dyn}$ one million times with values of $w$ sampled uniformly between its bounds, the worst-case cost obtained is $8.97 \times 10^{-6}$, which is within the solver tolerance of $\gamma^{\, *}$.  

\vspace{-4pt}
\subsection{Parameter Estimation}
\vspace{-4pt}

Finally, we solve the example of Section \ref{sec:ex3} with $n = 5$ time steps, where the goal is to estimate the tightest possible bounds $\overline{m}$ and $\underline{m}$ that could have resulted in the measurements
\begin{subequations}
    \begin{equation}
        \mbf{y} = \begin{bmatrix} -0.1 & 0.0 & 0.9 & 3.0 & 6.1 & 10.2 \end{bmatrix}^\top
    \end{equation}
under the control inputs
    \begin{equation}
        \mbf{u} = \begin{bmatrix} 1.0 & 1.0 & 1.0 & 1.0 & 1.0 \end{bmatrix}^\top 
    \end{equation}
\end{subequations}
It is also known that $x_{2,0}$ is $0$, but $x_{1,0}$ can only be inferred from $\mbf{y}$, $\mbf{u}$, and knowledge of the system dynamics.

The methodology of Section \ref{sec:local_red}, when applied to this problem, results in two additions to the optimal scenario set. These correspond to the two values of $\mbf{w}$ that lead to the respective worst-case solutions $\overline{m}^{*} = 1.01$ and $\underline{m}^* = 0.943$. The corresponding bound on the cost for the fit is $\gamma^{\, *} = 4.46 \times 10^{-3}$.

\section{Conclusions}

This paper introduced a local reduction algorithm that is capable of handling a wide range of min-max problems with existence constraints and non-unique solutions. The proposed methodology can be used to tackle numerous relevant robust optimal control and estimation problems, as was done with the examples considered on obstacle avoidance, control under input saturation, and parameter estimation. Local optimization solvers were used to obtain solutions that avoided any constraint violations while achieving a good performance in the cost, eliminating the need for the global optimizers that are frequently used in treating SIPs.  

The work presented here was limited to the computation of open-loop control input trajectories that guarantee constraint satisfaction over the whole of the prediction horizon. The method could potentially be extended to MPC implementations that would enable the use of information obtained at later time steps to improve the control performance.

%\begin{ack}
%Place acknowledgments here.
%\end{ack}

\bibliography{ifacconf}             % bib file to produce the bibliography
                                                     % with bibtex (preferred)
                                                   
%\begin{thebibliography}{xx}  % you can also add the bibliography by hand

%\bibitem[Able(1956)]{Abl:56}
%B.C. Able.
%\newblock Nucleic acid content of microscope.
%\newblock \emph{Nature}, 135:\penalty0 7--9, 1956.

%\bibitem[Able et~al.(1954)Able, Tagg, and Rush]{AbTaRu:54}
%B.C. Able, R.A. Tagg, and M.~Rush.
%\newblock Enzyme-catalyzed cellular transanimations.
%\newblock In A.F. Round, editor, \emph{Advances in Enzymology}, volume~2, pages
%  125--247. Academic Press, New York, 3rd edition, 1954.

%\bibitem[Keohane(1958)]{Keo:58}
%R.~Keohane.
%\newblock \emph{Power and Interdependence: World Politics in Transitions}.
%\newblock Little, Brown \& Co., Boston, 1958.

%\bibitem[Powers(1985)]{Pow:85}
%T.~Powers.
%\newblock Is there a way out?
%\newblock \emph{Harpers}, pages 35--47, June 1985.

%\bibitem[Soukhanov(1992)]{Heritage:92}
%A.~H. Soukhanov, editor.
%\newblock \emph{{The American Heritage. Dictionary of the American Language}}.
%\newblock Houghton Mifflin Company, 1992.

%\end{thebibliography}

%\appendix
%\section{A summary of Latin grammar}    % Each appendix must have a short title.
%\section{Some Latin vocabulary}              % Sections and subsections are supported  
                                                                         % in the appendices.
\end{document}